\documentclass[12pt]{amsart}

\usepackage{stmaryrd}

\usepackage{amssymb,amsmath,amsfonts,latexsym,tikz}
\newtheorem{theorem}{Theorem}
\newtheorem{definition}[theorem]{Definition}
\newtheorem{proposition}[theorem]{Proposition}

\theoremstyle{definition}
\newtheorem{example}[theorem]{Example}
\newtheorem{remark}[theorem]{Remark}

\renewcommand{\thefootnote}{\fnsymbol{footnote}}

\def\Z{\ensuremath{\mathbb{Z}}}
\def\Q{\ensuremath{\mathbb{Q}}} 
\def\P{\ensuremath{\mathbb{P}}}

\def\F{\ensuremath{\mathbb{F}}}

\def\C{\ensuremath{\mathbb{C}}}
\def\A{\ensuremath{\mathbb{A}}}
\def\R{\ensuremath{\mathbb{R}}}

\def\fp{\ensuremath{\mathfrak{p}}}
\def\fm{\ensuremath{\mathfrak{m}}}
\def\V{\ensuremath{\mathcal{V_\R}}}

\def\int{\mathrm{int}}

\def\an{\ensuremath{\mathrm{an}}}

\def\Gal{\ensuremath{\mathrm{Gal}}}
\def\Spec{\ensuremath{\mathrm{Spec \,}}}
\def\Spf{\ensuremath{\mathrm{Spf \,}}}

\def\<{\ensuremath{\langle}}
\def\>{\ensuremath{\rangle}}

\def\Link{\ensuremath{\mathrm{Link}}}

\title[Topology of nonarchimedean analytic spaces]{Topology of nonarchimedean analytic spaces and relations to complex algebraic geometry}

\author{\textsc{Sam Payne}}
\address{Yale University
Mathematics Department\\
10 Hillhouse Ave\\
New Haven, CT 06511 \\ U.S.A.}
\email{sam.payne@yale.edu}
\thanks{Partially supported by NSF DMS-1068689 and NSF CAREER DMS-1149054.}
\date{}

\begin{document}
\renewcommand*{\thefootnote}{\arabic{footnote}}

\maketitle

\begin{abstract}
    This note surveys basic topological properties of nonarchimedean analytic spaces, in the sense of Berkovich, including the recent tameness results of Hrushovski and Loeser.  We also discuss interactions between the topology of nonarchimedean analytic spaces and classical algebraic geometry.
    \end{abstract}

\setcounter{tocdepth}{1}
\tableofcontents

\section{Introduction}

\subsection{Complex algebraic geometry.}  At its most basic, classical complex algebraic geometry studies the common zeros in $\C^n$ of a collection of polynomials in $\C[x_1, \ldots, x_n]$.  Such an algebraic set may have interesting topology, but is not pathological.  It can be triangulated and admits a deformation retract onto a finite simplicial complex.  Furthermore, it contains an everywhere dense open subset that is a complex manifold, and whose complement is an algebraic set of smaller dimension.  Proceeding inductively, every algebraic set in $\C^n$ decomposes as a finite union of complex manifolds, and many of the deepest and most fundamental results in complex algebraic geometry are proved using holomorphic functions and differential forms, Hodge theory, and Morse theory on these manifolds.

\subsection{Beyond the complex numbers.}  Modern algebraic geometers are equally interested in the common zeros in $K^n$ of a collection of polynomials in $K[x_1, \ldots, x_n]$, for fields $K$ other than the complex numbers.  For instance, the field of rational numbers $\Q$ is interesting for arithmetic purposes, while the field of formal Laurent series $\C(\!(t)\!)$ is used to study deformations of complex varieties.  Like $\C$, such fields have natural norms.  For a prime number $p$, the $p$-adic norm $| \ |_p$ is given by writing a nonzero rational number uniquely as $\frac{p^a r}{s}$, with $a$ an integer, and $p$, $r$, and $s$ relatively prime, and then setting
\[
\left |\frac{p^ar}{s}\right |_p = p^{-a}.
\]
The $t$-adic norm $| \ |_t$ is given similarly, by writing a formal Laurent series uniquely as $t^a$ times a formal power series with nonzero constant term, and then setting
\[
\left |t^a \sum a_i t^i \right|_t = e^{-a}.
\]
One can make sense of convergent power series with respect to these norms, and it is tempting to work naively in this context, with ``analytic" functions given locally by convergent power series.  Difficulties arise immediately, however, for essentially topological reasons, unless the field happens to be $\C$.  The pleasant properties of analytic functions in complex geometry depend essentially on $\C$ being an archimedean field.

\subsection{What is an archimedean field?} The \emph{archimedean axiom} says that, for any $x \in K^*$, there is a positive integer $n$ such that $|nx| > 1$.  An archimedean field is one in which this axiom holds, such as the real numbers and the complex numbers.  However, there are essentially no other examples.  The archimedean axiom is satisfied only by $\C$, with powers of the usual norm, and restrictions of these norms to subfields.  In particular, the only complete archimedean fields are $\R$ and $\C$.  

\subsection{A nonarchimedean field is any other complete normed field}   We are not talking about the snake house or a rare collection of exotic creatures.  Nonarchimedean fields are basically the whole zoo.  Examples include the completion $\Q_p$ of $\Q$ with respect to the $p$-adic norm, and the field of formal Laurent series $\C(\!(t)\!)$.  Also, every field is complete and hence nonarchimedean with respect to the trivial norm, given by $|x| = 1$ for $x \in K^*$.

The norm on a nonarchimedean field extends uniquely to its algebraic closure.  The algebraic closure may not be complete\footnote{This is not difficult to see in examples.  For instance, the algebraic closure of $\C(\!(t)\!)$ is the field of Puiseux series $\C\{\!\{t\}\!\} = \bigcup_n \C(\!(t^{1/n})\!)$.  The exponents appearing in a Puiseux series have denominators bounded above, but these bounds need not be uniform on a Cauchy sequence.  For instance, the sequence $x_n = \sum_{j=1}^n t^{j + \frac{1}{j}}$ is Cauchy, but has no limit in $\C\{\!\{t\}\!\}$.}, but the completion of a normed algebraically closed field is again algebraically closed \cite[Proposition~3.4.1.3]{BGR84}.  So the completion of the algebraic closure of a nonarchimedean field is both nonarchimedean and algebraically closed.  A typical example is $\C_p$, the completion of the algebraic closure of $\Q_p$.

\subsection{The ultrametric inequality.} In a nonarchimedean field, a much stronger version of the triangle inequality holds.  The \emph{ultrametric inequality} says that 
\[
|x + y| \leq \max \{|x|, |y|\}, \mbox{ with equality if } |x| \neq |y|.
\]  
This property deserves a few moments of contemplation.  It implies that, if $y$ is a point in the open ball 
\[
B(x,r) = \{ y \in K \ | \ |y-x| < r\},
\]
then $B(x,r) = B(y,r)$.  In other words, every point in a nonarchimedean ball is a center of the ball.

\begin{remark}
Some authors define a nonarchimedean field to be any normed field that satisfies the ultrametric inequality and do not require that the field be complete with respect to the norm.  Since completeness is essential for analytic geometry, we maintain this additional hypothesis throughout.  
\end{remark}

\subsection{Nonarchimedean fields are totally disconnected.}  Because of the ultrametric inequality, the open ball $B(x,r)$ in a nonarchimedean field is closed in the metric topology.  Since these sets form a basis for the topology, the field $K$ is totally disconnected. Doing naive analysis on such a totally disconnected set is unreasonable.  For instance, if $f$ and $g$ are any two polynomials, then the piecewise defined function
\[
\Phi(x) = \left \{ \begin{array}{ll} f(x) & \mbox{ if $x \in B(0,1)$;} \\
								g(x) & \mbox{ otherwise,} \end{array} \right.	
\]
is continuous.  Even worse, this function $\Phi$ is ``analytic" in the naive sense that it is given by a convergent power series in a neighborhood of every point.

\subsection{Grothendieck topologies.} Let $K$ be a nonarchimedean field.  The affine space $K^n$ is totally disconnected in its metric topology, so a purely naive approach to analytic geometry over $K$ is doomed to fail.  One kludge is to discard the naive notion of topology.  

Let us return for a moment to the space of rational numbers $\Q$, which is totally disconnected in its metric topology.  The interval $[0,1]$ in $\Q$ is totally disconnected and non compact.  However, this totally disconnected set cannot be decomposed into a disjoint union of two open segments with rational centers and rational endpoints.   This suggests that $[0,1] \cap \Q$ is in some sense connected with respect to finite covers by rational intervals (even though it does decompose as a disjoint union of two open sets, each of which is an infinite union of rational intervals).  This suggests that one should restrict to considering finite covers by rational intervals (or at least covers with a finite refinement) in order to do analysis on $\Q$.  An approach like this can work, once one gives up the idea that an arbitrary union of open sets should be open.  Naive topology involving open sets and covers by open sets is then replaced by a \emph{Grothendieck topology}, consisting of a collection of covers satisfying certain axioms that are satisfied by the usual open covers in topology.

\subsection{Rigid analytic geometry.} John Tate developed a satisfying and powerful theory of nonarchimedean analytic geometry, based on sheaves of analytic functions in a Grothendieck topology on $K^n$, when $K$ is algebraically closed.\footnote{When $K$ is not necessary algebraically closed, Tate's theory uses a Grothendieck topology on $\overline K^n / \Gal(\overline K | K)$.  Experts will note that this is the set of closed points in the scheme $\Spec K[x_1, \ldots, x_n]$.}  His theory with this Grothendieck topology is called \emph{rigid analytic geometry}.  The name ``rigid" contrasts these spaces from the na\"{i}ve totally disconnected analytic spaces, which Tate called ``wobbly."  The fundamental algebraic object in the theory, the ring of convergent power series on the unit disc, is called the \emph{Tate algebra}.  Algebraic properties of the Tate algebra, including the fact that it is noetherian, play an essential role in all forms of nonarchimedean analytic geometry, whether one works in the rigid setting or follows the approach of Berkovich.  See \cite{BGR84} for a comprehensive treatment of the foundations of rigid analytic geometry.

\subsection{Filling in gaps between points.}  As mentioned above, one can develop a version of analysis on $\Q$ by replacing the metric topology with a suitable Grothendieck topology.  Nevertheless, most mathematicians prefer to add in new points that ``fill in the gaps" between the rational numbers and do analysis on the real numbers instead.  Note the fundamental absurdity of this construction.  Although $\Q$ is dense in $\R$, it has measure zero.  Once we have filled in the gaps, we can more or less ignore the points in $\Q$ when we do analysis.  What is added is so much larger than what we started with.  

In the late 1980s and early 1990s, Vladimir Berkovich developed a new version of analytic geometry over nonarchimedean fields.  At the heart of his construction is a topological space that fills in the gaps between the points of $K^n$, producing a path connected, locally compact Hausdorff space that contains $K^n$ with its metric topology, and $K^n$ is dense if the norm is nontrivial and $K$ is algebraically closed.\footnote{If the norm is nontrivial but $K$ is not algebraically closed, then $K^n$ may not be dense.  However, Berkovich's analytification also contains $\overline K^n / \Gal(\overline K | K)$ with its natural topology induced by the metric on $\overline K$, and this subspace is dense.}  The underlying algebra and analysis in Berkovich's theory are essentially the same as in rigid analytic geometry, but the topological space is different.  The subject of this note is the topology of the spaces appearing in Berkovich's theory, recent results on the tameness of these spaces, and relations between topological invariants of these spaces and more classical notions in complex algebraic geometry.  The first four sections are written for a general audience, while the final section, on relations to complex algebraic geometry, assumes some familiarity with the cohomology of algebraic varieties and contains an example illustrating the failure of Lefschetz theorems with integer coefficients on nonarchimedean analytic spaces.

\section{Nonarchimedean analytification}  \label{sec:analytification}

Nonarchimedean analytification is a functor from algebraic varieties (or, more generally, separated schemes of finite type) over a nonarchimedean field to analytic spaces in the sense of Berkovich.  For simplicity, we focus on the case of an affine variety.  Analytifications of arbitrary varieties are obtained by a natural gluing procedure from analytifications of affine varieties, which can be described as spaces of seminorms on coordinate rings.

\subsection{Seminorms on coordinate rings.}  Let $K$ be a nonarchimedean field.  Consider polynomials $f_1, \ldots, f_r \in K[x_1, \ldots, x_n]$, and let $X$ be the space of solutions to the corresponding system of equations.\footnote{In other words, $X$ is the Zariski spectrum $\Spec K[x_1, \ldots, x_n]/(f_1, \ldots, f_r)$, a locally ringed space whose underlying topological space is the set of prime ideals $\fp$ in this quotient ring, with the Zariski topology.  Yoneda's Lemma identifies this space with the functor that associates to a $K$-algebra $S$ the set of solutions to $f_1, \ldots, f_r$ in $S^n$.  In particular, if $L | K$ is an extension field then $X(L)$ is the set of points $y \in L^n$ such that $f_i(y) = 0,$ for $1 \leq i \leq r$.}  If $x = (x_1, \ldots, x_n)$ is a point in $X(K)$, then there is an associated seminorm on the quotient ring
\[
K[X] = K[x_1, \ldots, x_n] / (f_1, \ldots, f_r).
\]
Here, a seminorm on a ring is simply a function $| \ |$ from $K[X]$ to $\R_{\geq 0}$ that satisfies most of the usual axioms of a norm on a field, specifically that
\[
|fg| = |f| \cdot |g|
\]
and
\[
|f + g| \leq |f| + |g|.
\]
The seminorm $| \ |_x$ associated to a point $x$ in $X(K)$ is given simply by
\[
|f|_x = |f(x)|.
\]
We will only consider seminorms with the additional property that the restriction to $K$ is the given norm.  The given norm on $K$ is nonarchimedean, and it follows that any seminorm on $K[X]$ extending this norm also satisfies the ultrametric inequality
\[
|f + g| \leq \max \{ |f|, |g| \}, \mbox{ with equality if } |f| \neq |g|.
\]
The one significant difference between norms on fields and seminorms on rings is that a seminorm may take the value zero on a nonzero element of the ring.

\subsection{Analytification of affine varieties.}

We now describe the analytification of the affine variety $X$ over $K$ in terms of seminorms on its coordinate ring.

\begin{definition}
The analytification $X^\an$ is the space of all seminorms on the ring $K[X]$ that extend the given norm on $K$.\footnote{A word of caution is in order.  The system of polynomials $f_1, \ldots, f_r$ may not have any solutions defined over $K$.  Nevertheless, if $K[X]$ is nonzero then $X^\an$ is not empty.  One way to see this is to observe that the system $f_1, \ldots, f_r$ has solutions over the algebraic closure $\overline K$.  Since $K$ is complete, its norm extends uniquely to $\overline K$, and hence a solution with coordinates in $\overline K$ also determines a point of $X^\an$.  More generally, if $L | K$ is an extension field with a norm that extends the given one on $K$, then any solution to $f_1, \ldots, f_r$ with coordinates in $L$ determines a point of $X^\an$.} \footnote{An analogous definition over the complex numbers with its archimedean norm recovers classical complex analytic spaces.  If $X$ is an affine variety over $\C$ then the associated complex analytic space is naturally identified with the space of seminorms on $\C[X]$ whose restriction to $\C$ is the usual archimedean norm.  The bijection takes a point $x \in X(\C)$ to the seminorm $|f|_x = |f(x)|$.}
\end{definition}

\noindent We write $x$ for a point of $X^\an$, when thinking geometrically, and $| \ |_x$ for the corresponding seminorm on $K[X]$.  The topology on $X^\an$ is the subspace topology for the natural inclusion
\[
X^\an \subset (\R_{\geq 0})^{K[X]}.
\]
This is the coarsest topology such that, for each $f \in K[X]$, the function on $X^\an$ given by $x \mapsto |f|_x$ is continuous.

\begin{theorem} [\cite{Berkovich90}]
The topological space $X^\an$ is Hausdorff, locally compact, and locally path connected. The induced topology on the subset $X(K)$ of points with coordinates in $K$ is the metric topology, and if $K$ is algebraically closed with nontrivial valuation, then this subset is dense.
\end{theorem}

\noindent In this sense, $X^\an$ is a reasonable topological space on which to do analysis that ``fills in the gaps" between the points in the totally disconnected set $X(K)$.

\medskip

The remainder of Section~\ref{sec:analytification} addresses some of the richer structures on nonarchimedean analytic spaces, beyond the underlying topological space.  The casual reader may safely skip ahead to the examples in Section~\ref{sec:examples}.

\subsection{Structure sheaf and morphisms.} As an analytic space, $X^\an$ comes with much more structure than just a topology.  It has a sheaf of analytic functions given locally near each point by limits of rational functions that are regular at that point, and analytic spaces are objects in a category whose arrows are continuous maps such that the pull back of an analytic function under an analytic map is analytic.  There are furthermore well-behaved notions of open and closed embeddings, as well as flat, smooth, proper, finite, and \'etale morphisms in the category of analytic spaces.  One example of an \'etale morphism is given by pulling back the structure sheaf to a topological covering space, so the fundamental group of the underlying topological space gives essential information regarding the \'etale homotopy type of the analytic space. The details may be found in \cite{Berkovich93}.  Here, we are content simply to study the topological space underlying $X^\an$.

\subsection{Projection to the scheme.}  For those familiar with scheme theory or Zariski spectra of rings, we explain how $X^\an$ relates to $\Spec K[X]$, the set of prime ideals in $K[X]$ with its Zariski topology.

Let $x$ be a point in $X^\an$.  The set of functions $f \in K[X]$ such that $|f|_x = 0$ is a prime ideal $\fp$.  The seminorm $| \ |_x$ factors through a norm on the residue field $\kappa_\fp$, the fraction field of the quotient $K[X] / \fp$, whose restriction to $K$ is the given one.  
\begin{definition}
For any extension field $L|K$, let $\V(L)$ be the space of all norms on $L$ that extend the given norm on $K$.
\end{definition}
The map taking a point in the analytification to the kernel of the corresponding seminorm gives a natural surjection
\[
X^\an \rightarrow \Spec K[X]
\]
whose fiber over a point $\fp$ is $\V(\kappa_\fp)$.  In particular, the analytification decomposes as a disjoint union
\begin{equation} \label{decomposition}
X^\an = \bigsqcup_{\fp \in X} \V(\kappa_\fp).
\end{equation}
Note that the topology on the scheme $\Spec K[X]$ is never Hausdorff unless $X$ has dimension zero.  Its nonclosed points are the nonmaximal prime ideals $\fp \subset K[X]$, and the closure of $\fp$ is the irreducible variety $\Spec K[X] / \fp$.  The process of analytification produces a Hausdorff space by replacing each nonclosed point $\fp$ with the space of norms $\V(\kappa_\fp)$.  Each closed point of $\Spec K[X]$ is a maximal ideal $\fm \subset K[X]$.  The residue field $\kappa_\fm$ at a maximal ideal is algebraic over $K$ \cite[Proposition~7.9]{AtiyahMacdonald69}, so the norm on $K$ extends uniquely.  Therefore, there is a single point of $X^\an$ over each closed point of $\Spec K[X]$.

\subsection{A quotient description of $X^\an$.}
The decomposition above shows that each point of $X^\an$ is associated to a point of $\Spec K[X]$ together with a norm on its residue field that extends the given one on $K$.  For some purposes, rather than keeping track of all of these residue fields, it makes more sense to consider points defined over arbitrary normed extensions $L |K$.  One can still recover the analytification $X^\an$ by taking a quotient by an appropriate equivalence relation, as described below.  The Zariski spectrum $\Spec K[X]$ has an analogous description, in terms of natural equivalence classes of points over extension fields of $K$.  The key difference here is the role played by norms.

A normed extension $L|K$ is a field $L$ together with a norm $| \ | : L \rightarrow \R_{\geq 0}$ that extends the given norm on $K$.  We consider triples consisting of an extension field of $K$, a norm that extends the given one, and a point of $X$ over this normed extension, and the equivalence relation generated by setting
\[
(L, | \ |, x) \sim (L', | \ |', x')
\]
whenever there is an embedding $L \subset L'$ such that the restriction of $| \ |'$ to $L$ is $| \ |$ and $x$ is identified with $x'$ by the induced inclusion $X(L) \subset X(L')$.

\begin{proposition}
The analytification $X^\an$ is the space of equivalence classes of points of $X$ over normed extensions of $K$:
\[
X^\an = \{(L, | \ |, x) \} / \sim.
\]
\end{proposition}

\noindent Much of the recent progress in understanding the topology of nonarchimedean analytic spaces has come through logic and model theory, and this description of $X^\an$ in terms of equivalence classes of points over normed extensions is closest in spirit to the \emph{spaces of stably dominated types} that appear prominently in this context.  Note, however, that model theorists typically consider seminorms into ordered groups of arbitrary rank, not only the real numbers.

\section{Examples: affine line, algebraic curves, and affine plane} \label{sec:examples}

If $X$ has dimension 0 then $X^\an$ is equal to $X$, and both have the discrete topology.  We now consider the first nontrivial cases of analytifications.

\subsection{Analytification of the line: trivial norm.}  The simplest example to consider is the affine line
\[
\A^1 = \Spec K[y],
\]
in the case where the norm on $K$ is trivial.  Let $x$ be a point in $(\A^1)^\an$.  If $|y|_x$ is greater than 1 and
\[
f = a_0 + a_1y + \cdots + a_d y^d
\]
is a polynomial of degree $d$, then $|f|_x = |y|_x^d$.  Therefore, $| \ |_x$ is uniquely determined by $|y|_x$, and the limit, as $|y|_x$ goes to 1, is the trivial norm $\eta$ on the function field $K(y)$.  This gives an embedded copy of $[1, \infty)$ in $(\A^1)^\an$.  

Now, suppose $| \ |_x$ is not trivial, and $|y|_x \leq 1$.  Then $| \ |_x$ is less than or equal to 1 on all of $K[y]$, and the set of $f$ such that $|f|_x < 1$ is a nonzero prime ideal.  Each such ideal is generated by a unique irreducible monic polynomial $g \in K[y]$.  Given such a $g$ and a real number $t < 1$, there is a unique seminorm on $K[y]$ such that $|g| = t$; it takes $g^a \cdot h$ to $at$, for $h$ relatively prime to $g$.  The limit of these seminorms as $t$ goes to 1 is again the trivial norm $\eta$, and the limit as $t$ goes to $0$ is the closed point corresponding to the maximal ideal $\fm_g$ generated by $g$.  This gives a rough picture of $(\A^1)^\an$ as a sort of tree, with an infinite stem consisting of seminorms on $K(y)$ that are greater than 1 on $y$, and infinitely many branches that end in leaves corresponding to the irreducible polynomials in $K[y]$.  Equivalently, the leaves correspond to closed points in the scheme $\A^1$ over $K$, or elements of $\overline K / \Gal(\overline K | K)$.
\begin{center}
\begin{figure}[h!] \label{A1}
\includegraphics{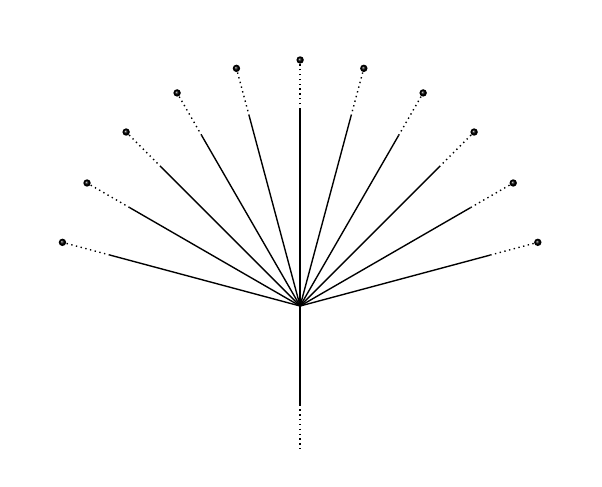}
\caption{The analytification of $ \A^1$ with respect to the trivial norm.}
\end{figure}
\end{center}
Some discussion of the topology on this tree is in order.  The topology on the subset where $|y| \leq 1$ is \emph{not} the cone over the discrete set $\overline K/ \Gal(\overline K | K)$.  Rather, it is an inverse limit of cones over finite subsets of $\overline K / \Gal(\overline K | K)$, so any neighborhood of $\eta$ contains all but finitely many of the branches.  To see this, note that for any $f \in K[y]$, the induced function $(\A^1)^\an \rightarrow \R \cup \infty$ taking $x$ to $|f|_x$ is 1 on almost all of the rays, all except those corresponding to irreducible factors of $f$.  Therefore, the preimage of any neighborhood of 1 in $\R_{\geq 0}$ contains all but finitely many of the branches, and these form a basis for the neighborhoods of $\eta$.

In this way, not only does $\V(K(y))$ fill in the gaps to connect the set of closed points in $\A^1_K$ with the discrete (metric) topology, it also interpolates between the metric topology and the cofinite (Zariski) topology in a subtle way.

\subsection{Analytification of the line: nontrivial norms.}

The analytification of $\A^1$ in the case of a nontrivial norm is again a tree, but now the set of branch points is dense.  At each branch point, the local topology is like the topology at $\eta$ in the analytification of the line with respect to the trivial norm.  It is described beautifully, and in detail, in Section~1 of \cite{Baker08b}.

\subsection{Analytification of curves.}
The analytification of an arbitrary smooth curve $X$ looks locally similar to that of the line.  If the norm is trivial, then $X^\an$ has finitely many open branches, corresponding to the points of the smooth projective model that are not in $X$, and the rest is an inverse limit of cones over finite subsets of $X(\overline K) / \Gal(\overline K | K)$.

If the norm on $K$ is nontrivial, then $X^\an$ is locally homeomorphic to $(\A^1)^\an$, but may have nontrivial global topology, as in the following example.  

\begin{center}
\begin{figure}[h!]
\includegraphics{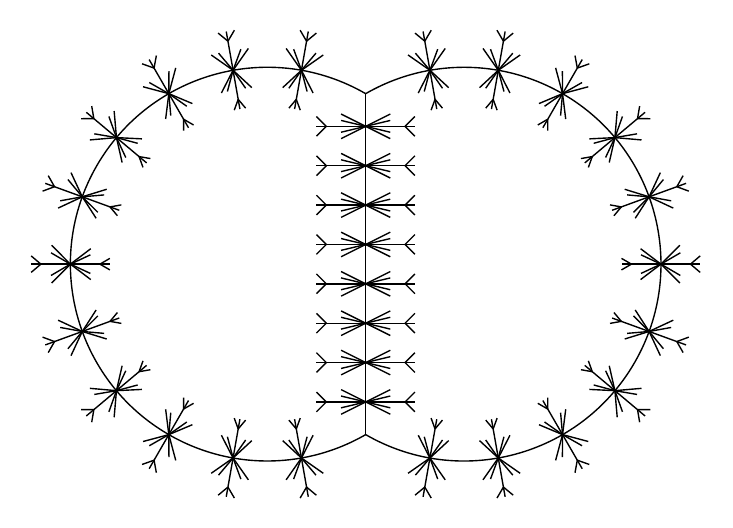}
\caption{The analytification of a smooth curve with respect to a nontrivial norm.}
\end{figure}

\end{center}

\noindent The dual graph of the special fiber of a semistable formal model embeds in $X^\an$ as a strong deformation retract.  For instance, an elliptic curve with bad reduction has a semistable formal model whose special fiber is a loop of copies of $\P^1$, and its analytification deformation retracts onto a circle.  Every finite graph occurs in this way, as the dual graph of the special fiber of a formal model, and hence as a deformation retraction of an analytic curve.

See \cite[Section~5]{BPR11} for further details on the structure theory of nonarchimedean analytic curves in the case where $K$ is algebraically closed.  The general case is similar; if $K$ is not algebraically closed then $X^\an$ is the analytification of the base change to the completion of the algebraic closure, modulo the action of $\Gal(\overline K | K)$.

In the case where $K$ has a countable dense subset, as is the case for $\Q_p$, $\C_p$, and $\widehat{\overline{\F_p}(\!(t)\!)}$, the topology of analytic curves over $K$ is locally modeled on that of the ``universal dendrite," each such curve admits a deformation retract onto a finite graph.  If the graph is planar then the curve admits an embedding in the euclidean plane $\R^2$, and if the graph is not planar then the curve embeds in $\R^3$ \cite{HrushovskiLoeserPoonen12}.

\subsection{Toward the analytification of the affine plane.}
Let us try to form a mental image of the analytification of the affine plane, using the discussion of curves, above, and the decomposition (\ref{decomposition}), in the case where the norm is trivial.  For another approach to visualizing the local topology of $(\A^2)^\an$, with illustrations, see \cite[Section 6.7]{Jonsson12}.

To start, note that $(\A^2)^\an$ contains the analytification of any plane curve, and the complement of the union of these analytic curves is the space $\V(K(x_1,x_2))$ of real norms on the function field $K(x_1,x_2)$ that are trivial on $K$.  Just as the closed points of a curve $X$ lie at the ends of infinite branches of $\V(K(X))$, the analytifications of curves in $\A^2$ lie in some sense at infinity, as limits of two-dimensional membranes in $\V(K(x_1,x_2))$.  Of course, the situation is somewhat more complicated, since the analytifications of distinct curves are joined at their points of intersection.

So, imagine a network of analytic curves at infinity, one for each curve in the plane, glued along leaves of the infinite branches corresponding to their points of intersection.  We now begin to describe how $\V(K(x_1,x_2))$ fills in the interior of this network.  Suppose $Y$ and $Z$ are curves in $\A^2$ that meet transversely at a point $x$.  Let $f$ and $g$ be defining equations for $Y$ and $Z$, respectively.  These can be interpreted as local coordinates on $\A^2$ at $x$, so any function $h \in K[x_1,x_2]$ can be expanded locally near $x$ as a power series in $f$ and $g$.

There is then a cone of ``monomial norms" in these local coordinates, for which the norm of a function depends only on the monomials in $f$ and $g$ that appear in its power series expansion, and these norms are defined as follows.  Let $v = (v_1, v_2)$ be a point in the cone $\R_{\geq 0 }^2$, and let $h$ be a function whose local power series expansion near $x$ is
\[
h = \sum a_{ij} f^i g^j.
\]
Then the norm corresponding to $v$ takes $h$ to
\[
|h|_v = \max_{i,j} \, \{ e^{-iv_1 - j v_2} \},
\]
where the maximum is taken over pairs $(i,j)$ such that $a_{ij} \neq 0$.  

The closure of this cone in $(\A^2)^\an$ is a copy of $(\R_{\geq 0} \cup \infty)^2$, joining the trivial norm $\eta$ on $K(x_1,x_2)$ to the trivial norms $\eta_Y$ and $\eta_{Z}$ on $K(Y)$ and $K(Z)$, respectively, and the point $x$.

\begin{center}
\begin{figure}[h!] \label{A2}
\includegraphics{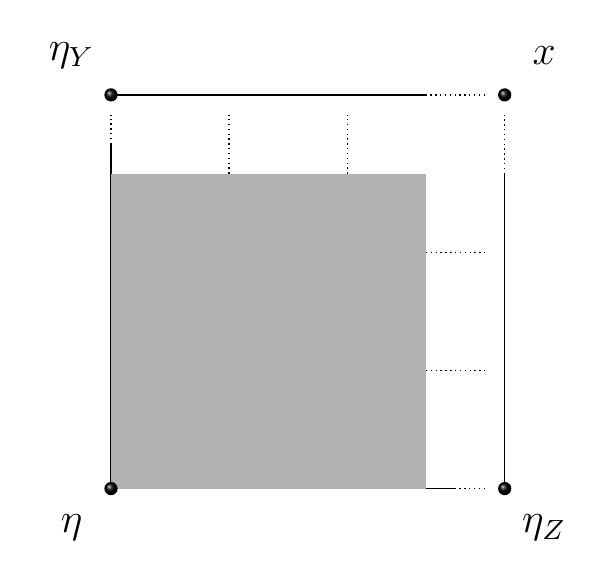}
\caption{The closure of a cone of monomial norms in $(\A^2)^\an$.}
\end{figure}
\end{center}

\noindent In the geometry of this cone, the limit of any ray with positive finite slope is $x$, so it is perhaps best imagined as a curved membrane stretching an infinite distance toward $x$, from the frame formed by the rays joining $\eta_Y$ and $\eta_{Z}$.

Understanding how all of these membranes fit together in $(\A^2)^\an$ is challenging, especially as one must also keep track of the topology in a neighborhood of $\eta$.  Moreover, we are still far from a full description of the underlying set of $(\A^2)^\an$.  All of the norms on $K(x_1,x_2)$ that are monomial in some system of local coordinates are of the simplest flavor in the sense of pure normed field theory; they satisfy Abhyankar's inequality\footnote{Recall that Abhyankar's inequality says that transcendence degree of the residue field extension plus the rank of the group extension given by the images of the norms is less than or equal to the transcendence degree of the total extension \cite{Abhyankar56}.} with equality.  There are many norms on $K(x_1,x_2)$ that are not monomial in any system of coordinates.  For instance, a pair of formal power series in $K\llbracket t \rrbracket$ define a formal germ of a curve in $\A^2$, and there is a seminorm on $K[x_1,x_2]$ obtained by pulling functions back to this germ and exponentiating minus the order of vanishing at $t = 0$.  If these power series are algebraically independent over $K$, then the image of this germ is Zariski dense in $\A^2$, and Abhyankar's inequality is strict. These points are outside of the infinite union of membranes described above.

Each seminorm corresponding to a point in $(\A^2)^\an$, including those where Abhyankar's inequality is strict, may be obtained as a limit of monomial seminorms in various systems of local coordinates (or even as a limit of seminorms corresponding to closed points), and the same is true in higher dimensions and on singular spaces, but the precise way that all of the pieces fit together becomes more and more difficult to describe.

\section{Tameness of analytifications}

Beyond the case of curves, which can be treated more or less by hand, it is not obvious that analytifications of algebraic varieties are not pathological topological spaces.  In his deep work on skeletons of formal models \cite{Berkovich99, Berkovich04}, Berkovich proved that analytifications of smooth varieties are locally contractible and have the homotopy type of a finite simplicial complex, provided that the norm is nontrivial.  However, the corresponding statements for singular varieties in positive and mixed characteristic, and for varieties over trivially normed fields, have so far eluded proof by such methods.  For instance, while it has been known for some time that the analytification of a smooth variety over a trivially normed field is contractible, the only known proof that they are locally contractible is quite recent and passes through model theory and spaces of stably dominated types.  This local contractibility is just one of the fundamental consequences of the tameness theorem of Hrushovski and Loeser \cite{HrushovskiLoeser10}.

\subsection{Analytic domains}

We now move toward describing the local structure of analytic spaces.  In algebraic geometry, affine open subvarieties are the basic building blocks.  An arbitrary algebraic variety is constructed by gluing its affine open subvarieties along open immersions, and each affine variety is determined by its ring of global regular functions.  Affinoid analytic domains play a similar role as basic building blocks in nonarchimedean analytic geometry.  An arbitrary nonarchimean analytic variety is constructed by gluing its affinoid analytic domains along immersions, and each affinoid domain is determined by its coordinate ring of global analytic functions.\footnote{Affine varieties and affinoid domains also share important cohomological and sheaf theoretic properties.  Each coherent sheaf on an affinoid domain is associated to a finitely generated module over its coordinate ring and has vanishing higher cohomology \cite[8.2.1 and 9.4.3]{BGR84}.  For the corresponding properties of affine varieties, see \cite[II.5.1 and III.3.5]{Hartshorne77}.}

We are particularly interested in analytifications of algebraic varieties, and now describe the affinoid domains in these spaces.  Let $X$ be an affine variety.  A typical affinoid domain in $X^\an$ can be realized by choosing a closed embedding $\iota: X \hookrightarrow \A^n$ and intersecting with the unit ball.  More precisely, if 
\[
\iota^* : K[y_1, \ldots, y_n] \rightarrow K[X]
\]
is the corresponding surjection of $K$-algebras, then the subset
\[
U = \{ x \in X^\an \ | \ |\iota^*(y_i)|_x \leq 1, \mbox{ for } 1 \leq i \leq n \}
\]
is an affinoid analytic domain.  As the terminology suggests, the affinoid analytic domain $U$ inherits the structure of a nonarchimedean analytic space from its inclusion in $X^\an$. 

The role of more general, not necessarily affine, open subvarieties in algebraic geometry is played by compact analytic domains in nonarchimedean analytic geometry.  The compact analytic domains are exactly the finite unions of affinoid analytic domains.  Every point in an analytic space has a basis of neighborhoods consisting of compact analytic domains, so understanding the topology of compact analytic domains is essential to understanding the local topology of analytic spaces.

\subsection{Skeletons of formal models}

Assume the norm on $K$ is nontrivial.  Then any compact analytic domain $U \subset X$ has formal models, which means that it can be realized as the ``analytic generic fiber" of a formal scheme.   Roughly speaking, a formal model of $U$ is given by a compatible system of local coordinates on $X$ such that the points in $U$ are exactly those whose local coordinates are in the valuation ring $R \subset K$.  In these local coordinates, one can then mod out by the maximal ideal $\fm \subset R$ to get the special fiber of the formal model, which is a scheme over the residue field $k = R/\fm$.

\begin{example}
Fix an affine embedding $X \subset \Spec K[x_1, \ldots, x_n]$.  Then we can ``clear denominators" in the coordinate ring $K[x_1, \ldots, x_n] / I_X$ to get a finitely presented $R$-algebra
\[
A = R[x_1, \ldots, x_n]/ (I_X \cap R[x_1, \ldots, x_n]).
\]
The scheme $\Spec A$ is an integral model of $X$.  Its generic fiber is $X$ and its special fiber is $\Spec A \otimes_R k$.  For any nonzero $\varpi \in \fm$, the $\varpi$-adic completion of this integral model is a formal model $\Spf \widehat A$ with analytic generic fiber
\[
U = \{ (x_1, \ldots, x_n) \in X^\an \ | \ |x_i| \leq 1 \mbox{ for } 1 \leq i \leq n \}.
\]
The special fiber of $\Spf \widehat A$ is the same as the special fiber of $\Spec A$.
\end{example}
  
\begin{remark}
Just as one can clear denominators on coordinate rings to get formal models of compact analytic domains, one can also clear denominators on morphisms to get morphisms between such formal models.  Raynaud famously proved that the category of quasi-compact and quasi-separated rigid analytic spaces is naturally identified with a localization of the category of quasi-compact formal schemes, in which admissible formal blowups (modifications of formal schemes that affect only the special fiber) are inverted.  See also \cite{Raynaud74b, BoschLutkebohmert93, BoschLutkebohmert93b, BLR95, BLR95b} for further details on formal schemes and their relation to analytic geometry. 
\end{remark}

If a compact analytic domain $U$ has a formal model with an especially nice special fiber\footnote{In residue characteristic zero, formal models with nice special fibers are constructed using the semistable reduction theorem \cite{KKMS}, a version of resolution of singularities for a one-parameter family of varieties.  Berkovich also proved local contractibility of smooth analytic spaces over nontrivially valued fields with positive residue characteristic using formal models with nice special fibers \cite{Berkovich99}, which he constructed via de Jong's theorem on alterations \cite{deJong96}.} then its topology is controlled by the combinatorics of the special fiber.

\begin{definition}
A formal model is strictly semistable if its special fiber is a reduced union of smooth varieties meeting transversally.  
\end{definition}

\noindent The dual complex of the special fiber of a strictly semistable formal model of $U$ is a regular $\Delta$-complex with one vertex for each irreducible component, one edge for each irreducible component of a pairwise intersection, one 2-face for each irreducible component of a triple intersection, and so on.  It has a canonical realization as a closed subset of $U$.  Roughly speaking, the dual complex consists of seminorms that are monomial in the local coordinates defining the model.  

\begin{theorem}[Berkovich]
If a compact analytic domain has a strictly semistable formal model, then it admits a strong deformation retract onto the dual complex of its special fiber.  In particular, it has the homotopy type of a finite simplicial complex.
\end{theorem}

\noindent Berkovich constructed skeletons and deformation retractions much more generally, for analytic spaces with polystable formal models, and used these skeletons to prove tameness results, including local contractibility, provided that such formal models exist, as is the case in residue characteristic zero \cite{Berkovich99, Berkovich04}.  See also \cite[Section~3 and Proposition~4.4]{Nicaise11} and \cite[Section~3]{MustataNicaise12} for an accessible treatment of the semistable case.  Other recent work of Nicaise and Xu shows that similar results hold for algebraic formal models whose special fibers have normal crossings but are not necessarily reduced \cite{NicaiseXu13}.

\begin{remark}
This discussion of skeletons of formal models assumes that the norm is nontrivial.  See \cite{Thuillier07} for a closely related construction of skeletons and deformation retractions associated to toroidal embeddings in the case where the norm is trivial.
\end{remark}

\subsection{Semialgebraic sets and tameness}

As mentioned above, Ber\-kovich proved that skeletons of sufficiently nice formal models control the topology of analytic spaces when the norm is nontrivial and such formal models exist, as is the case when the residue field has characteristic zero.  The recent work of Hrushovski and Loeser proves similar tameness results with no hypothesis on the normed field, by very different methods.  To state their tameness theorem, it is most helpful to talk about subsets of analytifications that are more general than compact analytic domains.

\begin{definition}
Let $X$ be an affine algebraic variety over $K$.  A semialgebraic subset $U \subset X^\an$ is a finite boolean combination of subsets of the form
\[
\{x \in X^\an \ | \ |f|_x \bowtie |g|_x^\lambda \},
\]
with $f, g \in K[X]$, $\lambda \in \R$, and $\bowtie \in \{ \leq, \geq, <, > \}$.
\end{definition}

\noindent  A semialgebraic set is \emph{definable} if the conditions defining the subset can be chosen such that some power of $\lambda$ is in $|K|$. 

By construction, every point in $X^\an$ has a basis of neighborhoods that are semialgebraic sets, and if the norm is nontrivial then these sets can be chosen to be definable.  The Gerritzen-Grauert theorem on locally closed immersions of affinoid varieties \cite[Theorem~7.3.5.1]{BGR84} guarantees that affinoid domains and, more generally, compact analytic domains in $X^\an$ are semialgebraic subsets.

We now state a basic version of the main result from \cite{HrushovskiLoeser10}.

\begin{theorem}
Let $U \subset X^\an$ be a semialgebraic subset.  Then there is a finite simplicial complex $\Delta \subset U$, of dimension less than or equal to $\dim(X)$, and a strong deformation retraction $U \times [0,1] \rightarrow \Delta$.
\end{theorem}

\noindent Since $\Delta$ is locally contractible, and the topology on $X^\an$ has a semialgebraic basis, it follows that $X^\an$ is locally contractible.  The homotopy type of the complex $\Delta$ is a fundamental invariant of a semialgebraic space $U$,  and many applications to complex geometry involve understanding the cohomology and fundamental groups of these complexes.

The approach of Hrushovski and Loeser does not involve the construction of nice models or toroidal compactifications.  It thereby avoids resolution of singularities and is insensitive to the residue characteristic.  As mentioned above, the proof involves a detailed study of spaces of stably dominated types, a notion coming from model theory \cite{HaskellHrushovskiMacpherson08}.  The case of curves is treated by hand, and the general case is a subtle induction on dimension, which involves birationally fibering an $n$-dimensional variety by curves over a base of dimension $n-1$.  In particular, the proof of the tameness theorem for a single variety in dimension $n$ requires a tameness statement for families of varieties in lower dimensions.  See the Bourbaki notes of Ducros \cite{Ducros13} for an excellent introduction to this work, and the original paper \cite{HrushovskiLoeser10} for further details.

\subsection{Limits of skeletons}

The finite simplicial complexes constructed in \cite{HrushovskiLoeser10}, which live inside semialgebraic sets as strong deformation retracts, are also called skeletons.  Hrushovski and Loeser prove much more than the existence of a single skeleton.  Each semialgebraic set $U$ of positive dimension contains infinitely many skeletons $\Delta_i$, with natural projections between them, and the semialgebraic set is recovered as the limit of this inverse system
\[
\varprojlim \Delta_i = U.
\]
Furthermore, each of these projections has a natural section, and the union $\varinjlim \Delta_i$ is the subset of $U$ consisting of points corresponding to Abhyankar seminorms.

Similar constructions involving limits of skeletons of formal models were considered earlier by Berkovich and in \cite{KontsevichTschinkel02, KontsevichSoibelman06}.  See \cite[Corollary~5.56 and Theorem~5.57]{BPR11} for an explicit treatment of such limits for curves.

\subsection{Limits of tropicalizations}

The topological space $X^\an$ can also be realized as a limit of finite polyhedral complexes using tropical geometry \cite{analytification, limits}.  In this approach, the polyhedral complexes are tropicalizations of algebraic embeddings of $X$ in toric varieties, with projections induced by toric morphisms that commute with the embeddings.  The construction of this inverse system is essentially elementary, at least in the quasiprojective case, but it does not lead to a proof of tameness.  It is not even known whether there exists a single tropicalization such that the projection from $X^\an$ is a homotopy equivalence, and the relation between these tropical inverse systems and the skeletons of Hrushovski and Loeser remains unclear.

\section{Relations to complex algebraic geometry}

The link between algebraic and analytic geometry over nonarchimedean fields is as close as the link between complex algebraic and complex analytic geometry.  Coherent algebraic sheaves have coherent analytifications, analytifications of \'etale algebraic morphisms are \'etale, and there are comparison theorems for $\ell$-adic \'etale cohomology \cite{Berkovich90, Berkovich93}.  It is not at all surprising, then, that nonarchimedean analytic techniques are powerful for studying algebraic varieties over nonarchimedean fields.

However, nonarchimedean analytic techniques are also powerful for studying algebraic varieties over the complex numbers.  One reason is simple: nonarchimedean fields such as $\C_p$ and the completion of $\C\{\!\{t\}\!\}$ are isomorphic to $\C$ as abstract fields.  The isomorphism is not explicit or geometric, but elimination of quantifiers for algebraically closed fields implies that any two uncountable algebraically closed fields of the same cardinality and characteristic are isomorphic \cite[Proposition~2.2.5]{Marker02}.  In particular, whenever one can use nonarchimedean analytic techniques to produce a variety over $\C_p$ with a certain collection of algebraic properties, it follows that there exists a variety over $\C$ with the same collection of properties.

Perhaps more surprisingly, one can also get significant mileage by studying analytifications of open and singular complex varieties with respect to the trivial norm and their semialgebraic subsets.  See the discussion of Milnor fibers and analytic links, below.

\subsection{Tropical linear series} Many applications of nonarchimedean analytic spaces in complex geometry involve less information than the full structure sheaf, but more than the mere topological space.  Tropical geometry resides firmly in this intermediate realm.  For instance, if $X$ is a curve then $\V(K(X))$, the complement in $X^\an$ of the set of closed points of $X$, inherits a natural metric.  Through the tropical Riemann-Roch theorem \cite{BakerNorine07, GathmannKerber08, MikhalkinZharkov08}, Baker's specialization lemma and its generalizations \cite{Baker08, AminiBaker12, AminiCaporaso13}, the nonarchimedean Poincar\'e-Lelong formula \cite{ThuillierThesis, BPR11}, and the theory of harmonic morphisms of metric graphs \cite{BakerNorine09, ABBR13}, this metric is a powerful tool in the study of linear series on algebraic curves. It has been used to characterize dual graphs of special fibers of regular semistable models of curves of a given gonality \cite{Caporaso12}, to compute the gonality of curves that are generic with respect to their Newton polygon \cite{CastryckCools12}, to characterize the Newton polygons of Brill-Noether general curves in toric surfaces \cite{Smith14}, to bound the gonality of Drinfeld modular curves \cite{CornelissenKatoKool12}, and to give new proofs of the Brill-Noether and Gieseker-Petri theorems \cite{tropicalBN, tropicalGP}.

In the remaining sections, we survey some of the applications of nonarchimedean analytic geometry to classical complex varieties that involve only the topology of analytifications.

\subsection{Singular cohomology}

Let $X$ be an algebraic variety over the complex numbers, and let $X(\C)$ be the associated complex analytic space.  Recall that Deligne defined a canonical mixed Hodge structure on the rational cohomology $H^*(X,\Q)$, and one part of this structure is the weight filtration
\[
W_0 H^k(X(\C), \Q) \subset \cdots \subset W_{2k}H^k(X(\C), \Q) = H^k(X(\C),\Q),
\]
which is strictly functorial for algebraic morphisms.  This means that if $f: X' \rightarrow X$ is a morphism then
\[
f^*(H^k(X(\C), \Q)) \cap W_j H^k(X'(\C), \Q) = f^* (W_jH^k(X(\C),\Q)).
\]
If $X$ is smooth and compact then $W_{k-1}H^k(X,\Q) = 0$ and $W_kH^k(X, \Q) = H^k(X,\Q)$.  In other words, $H^k(X(\C), \Q)$ is of pure weight $k$.  In the general case, where $X$ may be singular and noncompact, the graded pieces of $H^k$ of weight less than $k$ encode information on the singularities of $X$, while the pieces of weight greater than $k$ encode information about the link of the boundary in a compactification.  This is perhaps best illustrated by an example.

\begin{example} \label{punctured}
Consider a curve $X$ of geometric genus 1, with three punctures and a single node $x$, as shown.  
\vspace{10 pt}
\begin{center}
\begin{figure}[h!]
\includegraphics{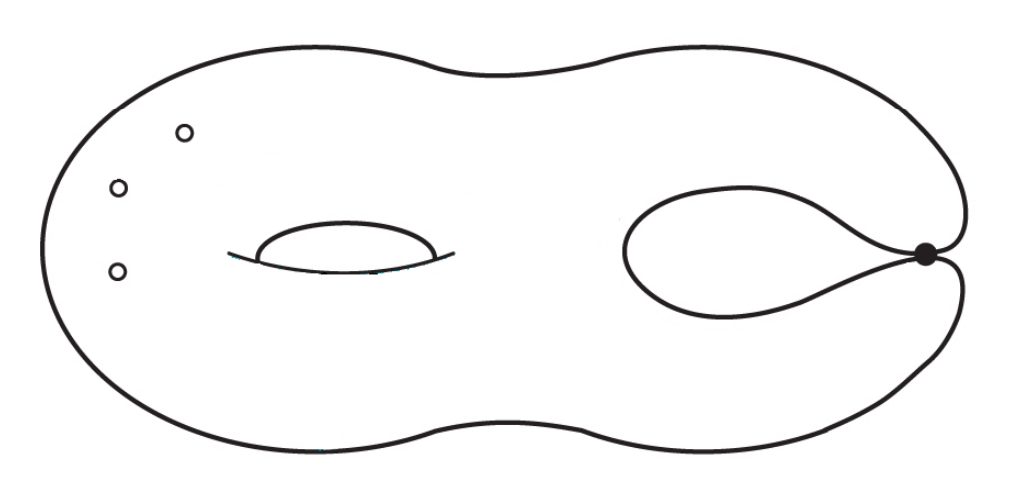}
\caption{A nodal curve with three punctures}
\end{figure}
\end{center}
Its normalization $\widetilde X$ is obtained by resolving the node.  The homology $H_1(X(\C),\Q)$ is generated by a loop through the node, two loops around the doughnut on the left, and two loops around punctures.  In the corresponding dual basis for $H^1(X(\C),\Q)$, the class dual to the loop through the node generates $W_0 H^1(X(\C),\Q)$.

We now consider the nonarchimedean analytification of $X$ with respect to the trivial norm on $\C$.  As in the examples from Section~\ref{sec:examples}, the analytification of the normalization $\widetilde X$ of $X$ is an infinite tree, with three unbounded branches corresponding to the punctures.  The remaining branches end in leaves, corresponding to the closed points of $\widetilde X$ and, in $X^\an$, two of these closed points are identified at the node $x$.  The analytification is therefore as shown here.

\begin{center}
\begin{figure}[h!] 
\includegraphics{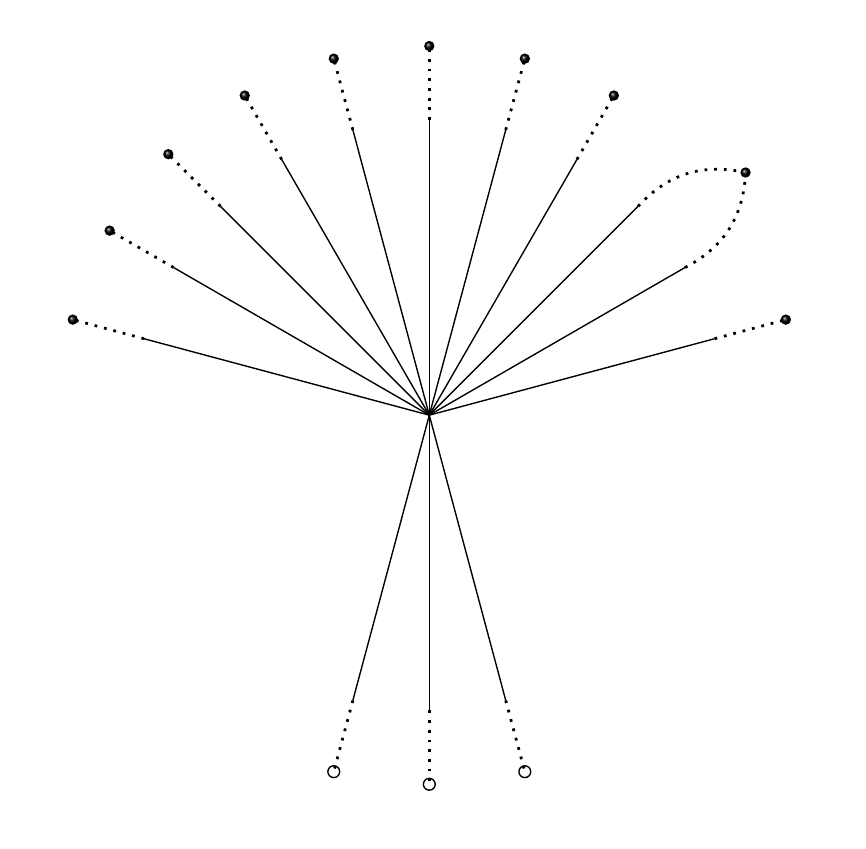}
\caption{The analytification of a nodal curve with three punctures. }
\label{nodal}
\end{figure}
\end{center}

\end{example}

Note that the identification of two closed points in $\widetilde X$ creates an extra loop in both $X(\C)$ and $X^\an$ and that, on $X(\C)$, the new class in $H^1$ has weight zero.

\begin{theorem}[\cite{Berkovich00}]
Let $X$ be a complex algebraic variety, and let $X^\an$ be its analytification with respect to the trivial norm on $\C$.  Then there is a natural isomorphism
\[
H^*(X^\an, \Q) \cong W_0 H^*(X(\C), \Q).
\]
\end{theorem}

A similar result holds for varieties defined over a local field, such as $\Q_p$.  In these cases, $H^k(X^\an,\Q_\ell)$ is canonically identified with the weight zero piece of the monodromy filtration on $H^k_{\acute{e}t}(X,\Q_\ell)$.

Singular cohomology of skeletons and their relations to weight filtrations have also appeared in the tropical geometry literature, for instance in \cite{Hacking08, HelmKatz12, KatzStapledon12}.  The key fact is that tropicalizations are skeletons in the case where all initial degenerations are smooth and irreducible, and there are natural parametrizing complexes for tropicalizations that are skeletons more generally, in the ``sch\"on" case, where all initial degenerations are smooth, but possibly reducible.  See \cite{Gubler13} for details on the relation between tropicalizations, initial degenerations, and formal models.

\subsection{Beyond rational cohomology}
There is far more information in the topology of $X^\an$ than just its rational cohomology.  For instance, a totally degenerate Enriques surface has no $\ell$-adic rational cohomology of weight zero in degree above zero, but its analytification is not contractible.  It has the homotopy type of $\R\P^2$.  In this case, the fundamental group of $X^\an$ agrees with the \'etale fundamental group of $X$.  The nature of the relationship between the fundamental group and higher homotopy groups of an analytification and the algebraic invariants of the variety is not known in general.  Subtleties appear already when one examines torsion in the cohomology of analytifications.  For instance, the naive generalization of the Lefschetz Hyperplane Theorem does not hold for cohomology of nonarchimedean analytic spaces, as explained in the next section.

\subsection{Lefschetz Hyperplane Theorems}
Suppose $D$ is a hyperplane section of a smooth projective algebraic variety over $K$.  The $\ell$-adic Lefschetz theorem \cite[Section 4.1.6]{Deligne80} says that the natural restriction maps
\[
H^i_{\acute{e}t}(X, \Z_\ell) \rightarrow H^i_{\acute{e}t}(D, \Z_\ell)
\]
are isomorphisms for $i < \dim D$ and injective for $i = \dim D$.  By tensoring with $\Q_\ell$ and applying the weight zero comparison theorems, the compatibility of restrictions with weight filtrations, and the universal coefficients theorem in singular cohomology, it follows that the natural maps
\[
H^i(X^\an,\Q) \rightarrow H^i(D^\an, \Q)
\]
are also isomorphisms for $i < \dim D$ and injective for $i = \dim D$.

Such Lefschetz Theorems also hold in the singular cohomology of complex varieties, and in this classical context they can be extended to integral cohomology and even homotopy groups.  Furthermore, the hyperplane section can be replaced by an arbitrary \emph{ample} divisor, i.e. a divisor $D$ such that $mD$ is a hyperplane section for some positive integer $m$.  The proof is by generalized Morse theory; a smooth complex variety of dimension $n$ can be obtained from any ample divisor by adding cells of dimension at least $n$ \cite{Bott59}.  Similar phenomena appear sometimes in tropical and nonarchimedean analytic geometry.  See, for instance, \cite{AdiprasitoBjorner14} for Lefschetz hyperplane theorems on locally matroidal tropical varieties.

\begin{example}
The analytification of the Jacobian of a totally degenerate curve of genus $g$ has a skeleton which is a real torus of dimension $g$.  This skeleton is naturally identified with the tropical Jacobian of the skeleton of the curve \cite{BakerRabinoff13}.  Furthermore, the image of the analytification of the theta divisor is the tropical theta divisor, and the tropical Jacobian is obtained from this tropical theta divisor by attaching a single cell of dimension $g$ \cite{MikhalkinZharkov08}.  This suggests that Lefschetz theorems may hold for the inclusion of the analytic theta divisor in the analytic Jacobian for integral cohomology and homotopy groups, though it seems to be unknown, in general, whether the projection from the analytification of the theta divisor to the tropical theta divisor is a homotopy equivalence.
\end{example}

\noindent However, it is not true that the inclusion of the analytification of an ample divisor in a variety of dimension $n$ induces isomorphisms on integral cohomology groups in degrees up to $n-2$, as the following example shows.

\begin{example}
Let $E$ be an elliptic curve over the complex numbers, with a $2$-torsion point $q$.  Let $E[2]$ be the $2$-torsion subgroup, which we view as a divisor of degree $4$ on $E$.  The abelian 3-fold $X = E \times E \times E$ has an ample divisor with simple normal crossings $D = p_1^*(E[2]) \cup p_2^*(E[2]) \cup p_3^* (E[2])$.  The dual complex $\Delta(D)$ has $12$ vertices, $48$ edges, and $64$ 2-faces; it is simply connected and homotopy equivalent to a wedge sum of $27$ spheres.  Note that $(q,q,q)$ acts on $X$ by a fixed point free involution that preserves $D$, and the induced action on $\Delta(D)$ is also fixed point free.  The quotient $(X',D')$ is again an abelian 3-fold with an ample divisor with simple normal crossings, and $\Delta(D) \rightarrow \Delta(D')$ is a universal cover.  In particular, 
\[
H^1(\Delta(D'), \Z) \, \cong \, \pi_1(\Delta(D')) \, \cong \, \Z/2\Z.
\]

Let $K$ be the completion of the algebraic closure of $\C(\!(t)\!)$.  We claim that there is a smooth ample divisor $H$ in $X'_K$, linearly equivalent to $D'$, whose analytification is homotopy equivalent to $\Delta(D')$.  To prove the claim, we first show that $D'$ is basepoint free.  Note that $D'$ is the pullback of a divisor $D''$ on $E/q \times E/q \times E/q$, and $D''$ is the union of the 3 pullbacks of a divisor of degree 2 on the elliptic curve $E/q$.  Any divisor of degree 2 on an elliptic curve is basepoint free.  It follows that $D''$ and $D'$ are basepoint free, since basepoint freeness is preserved by pullbacks and unions.  The total space of a general pencil containing $D'$ is smooth, by the Bertini Theorems, and completing this pencil produces a divisor $H$ in $X_K$ with a semistable model whose special fiber is $D'$.  In particular, $H$ is a smooth divisor in $X_K$, equivalent to $D'$, with skeleton $\Delta(D')$, as required.

The analytification of $X_K$ is contractible, because $X$ is smooth and defined over the subfield $\C$ on which the $t$-adic norm is trivial.  However, $H_1(H^\an, \Z)$ and $\pi_1(H^\an)$ are isomorphic to $\Z/2\Z$.  Thus, even though both $X$ and $H$ are smooth, and $H$ is ample in $X$, the natural maps
\[
H^1(X_K^\an, \Z) \rightarrow H^1(H^\an, \Z) \mbox{ \ \ and \ \ } \pi_1(H^\an) \rightarrow \pi_1(X_K^\an)
\]
are not isomorphisms.
\end{example}

\subsection{Specialization from analytic points to closed algebraic subsets}

Many semialgebraic subsets of analytifications that are of interest for complex algebraic geometry are defined in terms of analytic points that are near closed algebraic subsets.  The notion of \emph{specialization} captures the rough idea of a point of $X^\an$ being close to a point, or more generally a Zariski closed subset of $X$.

Let $x \in X^\an$ be a point and let $Z \subset X$ be a Zariski closed subset.  

\begin{definition}
We say that $x$ specializes into $Z$ if it is represented by a point $x \in X(K)$ for a valued extension $K | \C$ with valuation ring $R$ such that the inclusion $\iota_x : \Spec K \rightarrow X$ extends to a morphism
\[
\widetilde \iota: \Spec(R) \rightarrow X,
\]
which maps the closed point of $\Spec R$ into $Z$.
\end{definition}

\noindent If $X \subset \C^n$ is affine, then the set of points specializing into $Z$ is defined by the conditions $|x_i| \leq 1$ for $1 \leq i \leq n$ and $|f| < 1$ for $f \in I_Z$.  Given that the coordinate functions $x_i$ have norm less than or equal to 1, it suffices to check that $|f_j| < 1$ for some finite generating set $\{f_j\}$ of $I_Z$, so the set of points specializing into $Z$ is semialgebraic.

\begin{example}
If $x$ is a closed point of $\A^1$, then the set of points in the analytification of $\A^1$ specializing to $x$ is the open ray in Figure~\ref{A1} pointing from the central vertex toward $x$, together with $x$ itself. 
\end{example}

\begin{example}
If $x$ is the closed point of $\A^2$ shown in Figure~\ref{A2}, then the entire open cone of monomial valuations specializes into $x$.  The open rays pointing from $\eta$ to $\eta_X$ and $\eta_{X'}$ specialize into $X$ and $X'$, respectively, but not into $x$.
\end{example}

\begin{example}
If $x$ is the node in Figure~\ref{nodal} then both open rays pointing from the central vertex toward $x$ specialize to $x$, as does $x$ itself.
\end{example}

\subsection{The analytic Milnor fiber}
One prime example of a semialgebraic construction in nonarchimedean analytic geometry analogous to a topological construction in complex geometry is the \emph{analytic Milnor fiber} of Nicaise and Sebag.  Although their construction works more generally, we consider for simplicity the Milnor fiber of a single point in a hypersurface, and choose coordinates so that this point is the origin $0 \in \A^n$.

Let $X \subset \C^n$ be the vanishing locus of a polynomial $f$, and assume that $X$ contains $0$.  Recall that the classical Milnor fiber of $(X,0)$ is the intersection of the locus of points $x \in \C^n$ such that $f(x)$ has fixed argument with a Euclidean sphere of radius $\epsilon \ll 1$ centered at $0$.  See \cite{Milnor68} for further details.

We now define the nonarchimedean analytic Milnor fiber.  Note that the polynomial $f$ defines a morphism $f^\an$ from the analytification of $\A^n$ to the analytification of $\A^1 = \Spec \C[y]$.  Fix some $0 < \epsilon < 1$, and let $z$ be the unique point of $(\A^1)^\an$ such that $|y|_z = \epsilon$.

\begin{definition}
The analytic Milnor fiber of $(X,0)$ is the subset in the analytification of $\A^n$ over $\C(\!(t)\!)$ consisting of points $x$ that specialize to $0$ such that $f^\an(x) = z$.
\end{definition}

\noindent The condition that $x$ specializes to $0$ exactly means that $x$ lies in the open unit polydisc
\[
D = \{ (x_1, \ldots, x_n) \ | \ |x_i| < 1 \mbox{ \ for \ } 1 \leq i \leq n\}.
\]
The analytic Milnor fiber is semialgebraic and definable over the normed extension $\C(\!(t)\!) | \C$ in which $|t| = \epsilon$; after base change to $\C(\!(t)\!)$, it is just the closed analytic subvariety of $D$ defined by $f = t$.  Therefore, its base change to the algebraic closure of $\C(\!(t)\!)$ carries an action of the absolute Galois group, which induces an action on its $\ell$-adic \'etale cohomology.  Nicaise and Sebag show that the $\ell$-adic \'etale cohomology of the analytic Milnor fiber, with the action of the procyclic generator of this Galois group, is canonically identified with the $\ell$-adic singular cohomology of the classical Milnor fiber, with its monodromy action \cite{NicaiseSebag07}.

There are multiple advantages to the approach of Nicaise and Sebag.  First, their definition makes sense in much greater generality, for varieties over an arbitrary field with the trivial norm, such as $\overline \F_p$, where one does not have the complex topology to work with.  Furthermore, over $\C$, their construction puts the Milnor fiber and its monodromy action in a context (semialgebraic sets of smooth rigid varieties) where motivic integration makes sense \cite{LoeserSebag03}.  See also \cite{HrushovskiLoeser11}, which recasts the resulting interpretation of the motivic zeta function in the framework of Hrushovski and Kazhdan \cite{HrushovskiKazhdan06}.  All of this work opens up possibilities for a conceptual approach to the monodromy conjectures of Igusa \cite{Igusa75} and of Denef and Loeser \cite{DenefLoeser98}.

\subsection{The analytic link of a singularity}

Another such construction is the analytic link of a point $x \in X(\C)$.  The classical link of $x$ is a topological space $\Link(X(\C),x)$ obtained by embedding an affine neighborhood of $x$ in $\C^n$ and intersecting with a small sphere centered at $x$.  It can be given a piecewise linear structure, by triangulating $X(\C)$ with $x$ as a vertex and taking the link of this vertex in the resulting simplicial complex, and this triangulated space is well-defined up to piecewise linear homeomorphism.  Fundamental groups of links have been particularly prominent in recent research activity; see the survey article \cite{Kollar13}.

\begin{definition}
The analytic link $\Link(X^\an,x)$ is the semialgebraic set in $X^\an \smallsetminus x$ consisting of points $x'$ that specialize to $x$.
\end{definition}

\noindent Roughly speaking, this means that the points of the analytic link are those close to $x$, but not $x$ itself.  The analytic link of an arbitrary closed algebraic subset is defined similarly, and behaves like a deleted tubular neighborhood in classical complex geometry.

Advantages of the analytic link include the fact that its construction is canonical, not depending on any choice of local embedding or triangulation, and that it carries the additional structure of an analytic space. As for the Milnor fibers discussed above, comparison theorems should give canonical isomorphisms\footnote{The analogous comparison theorem for Milnor fibers was proved by Nicaise and Sebag \cite{NicaiseSebag07}, but to the best of our knowledge no such comparison for links has appeared in the literature.}
\[
H^*_{\acute{e}t}(\Link(X^\an,x), \Q_\ell) \simeq H^*(\Link(X(\C),x), \Q_\ell).
\]

The weight zero piece of the \'etale cohomology of the link corresponds to the singular cohomology of the underlying topological space, as in \cite{Berkovich00, Nicaise11}, so cohomological properties of singularities induce cohomological conditions on the topology of $\Link(X^\an, x)$.  For instance, if $(X,x)$ is an isolated rational singularity, then $\Link(X^\an,x)$ has the rational homology of a point, and if $(X,x)$ is an isolated Cohen-Macaulay singularity of dimension $n$, then $\Link(X^\an,x)$ has the rational homology of a wedge sum of spheres of dimension $n-1$.  These conditions can also be interpreted in terms of a log resolution; the dual complex of the exceptional divisor in a log resolution has the same homotopy type as the analytic link \cite{Thuillier07}.

Examples show that analytic links of rational singularities are not necessarily simply connected and may have torsion in their singular homology.  See \cite[Example~8.1]{boundarycx} and \cite{KapovichKollar11}.  However there are stronger natural conditions on singularities that do imply contractibility of the analytic link.  For instance, toric singularities and finite quotient singularities have contractible analytic links \cite{Stepanov06, KerzSaito12}.

De Fernex, Koll\'ar, and Xu recently proved the strongest results in this direction, showing that analytic links are contractible for isolated log terminal singularities.  Log terminal singularities are a subset of rational singularities that contain all toric and finite quotient singularities and appear naturally in birational geometry.  The proof in \cite{deFernexKollarXu13} involves the study of dual complexes of exceptional divisors not only for log resolutions, but for divisorial log terminal partial resolutions, and running a carefully chosen minimal model program, while keeping track of how these dual complexes transform at each step.

\subsection{The analytic link at infinity}

Analytic links of singularities are a special case of analytic links at infinity.  Again, we consider the analytification of a variety $X$ over the complex numbers with respect to the trivial norm but now, instead of studying the singularity at a point $x$, we examine the failure of $X$ to be compact.

\begin{definition}
The link at infinity $\Link_\infty(X^\an)$ is the semialgebraic subset of $X^\an$ consisting of points that do not specialize to any point of $X$.
\end{definition}

\noindent The link at infinity consists of points defined over valued extensions $K | \C$ that are not defined over the valuation ring $R \subset K$.  If $X$ is affine, then $x \in \Link_\infty(X^\an)$ if and only if $|f|_x > 1$ for some polynomial $f \in \C[X]$.  By the triangle inequality, it is enough to consider $f$ in some finite generating set, such as the coordinate functions for some embedding $X \subset \C^n$.

\begin{example}
Suppose $X$ is compact.  Then  $\Link_\infty(X^\an)$ is empty and, for any point $x \in X$, the space $\Link(X^\an, x)$ coincides with $\Link_\infty(X \smallsetminus x)$.
\end{example}

\begin{example}
Suppose $X$ is the nodal curve with three punctures discussed in Example~\ref{punctured}.  Then $\Link_\infty(X^\an)$ consists of the three open rays in Figure~\ref{nodal} pointing at the three punctures.
\end{example}

\begin{remark}
The link at infinity is also characterized as the link of the boundary in any compactification of $X$.  In other words, if $\overline X$ is a compactification of $X$ with boundary $\partial \overline X = \overline X \smallsetminus X$ the $\Link_\infty(X^\an)$ is the set of points of $X^\an$ that specialize to a point in $\partial \overline X$.
\end{remark}

\noindent The link at infinity is related to compactifications in much the same way that the link of an isolated singularity is related to resolutions.  Suppose $X$ is smooth and $\overline X$ is a smooth compactification whose boundary
\[
\partial \overline X = \overline X \smallsetminus X
\]
is a divisor with simple normal crossings.  Then Thuillier's constructions gives a canonical homotopy equivalence from $\Link_\infty(X)$ to the dual complex $\Delta(\partial \overline X)$ \cite{Thuillier07}.  One can then use excision exact sequences and Poincar\'e duality to identify the rational cohomology of $\Link_\infty(X^\an)$ with the top graded piece of the weight filtration on the cohomology of $X(\C)$,
\[
H^i(\Link_\infty(X^\an), \Q) \cong \mathrm{Gr}^W_{2\dim X} H^{2 \dim X - i}(X(\C), \Q).
\]
Such constructions with dual complexes of boundary divisors in smooth compactifications are surveyed in \cite{boundarycx}, but the main ideas were introduced and studied much earlier by Danilov \cite{Danilov75}.  This framework can be extended to toroidal compactifications of smooth Deligne-Mumford stacks, and this generalization has been applied to compute top weight cohomology groups for many moduli spaces of stable curves with marked points, including $\mathcal{M}_{1,n}$ for all $n$, using an interpretation of the dual complex of the boundary of the Deligne-Mumford compactification as a moduli space for stable tropical curves \cite{acp, cgp}.

\bigskip

\noindent \textbf{Acknowledgments.}  I thank P.~Achinger, M.~Baker, D.~Cart\-wright, L.~Fantini, S.~Grushevsky, J.~Huh, A.~Kontorovich, J.~Nicaise, D.~Ranganathan, D.~Speyer, W.~Veys, and the referee for helpful comments and important corrections to earlier versions of this article, and am grateful to D.~Ranganathan also for assistance in preparing the figures.

\bibliographystyle{amsalpha}
\bibliography{math}  

\providecommand{\bysame}{\leavevmode\hbox to3em{\hrulefill}\thinspace}
\providecommand{\MR}{\relax\ifhmode\unskip\space\fi MR }
\providecommand{\MRhref}[2]{%
  \href{http://www.ams.org/mathscinet-getitem?mr=#1}{#2}
}
\providecommand{\href}[2]{#2}
\begin{thebibliography}{KKMSD73}

\bibitem[AB12]{AminiBaker12}
O.~Amini and M.~Baker, \emph{Linear series on metrized complexes of algebraic
  curves}, preprint arXiv:1204.3508, 2012.

\bibitem[AB14]{AdiprasitoBjorner14}
K.~Adiprasito and A.~Bj{\"o}rner, \emph{Filtered geometric lattices and
  {L}efschetz section theorems over the tropical semiring}, preprint
  arXiv:1401.7301, 2014.

\bibitem[ABBR13]{ABBR13}
O.~Amini, M.~Baker, E.~Brugall\'e, and J.~Rabinoff, \emph{Lifting harmonic
  morphisms of tropical curves, metrized complexes, and {B}erkovich skeleta},
  preprint arXiv:1303.4812, 2013.

\bibitem[Abh56]{Abhyankar56}
S.~Abhyankar, \emph{On the valuations centered in a local domain}, Amer. J.
  Math. \textbf{78} (1956), 321--348.

\bibitem[AC13]{AminiCaporaso13}
O.~Amini and L.~Caporaso, \emph{Riemann--{R}och theory for weighted graphs and
  tropical curves}, Adv. Math. \textbf{240} (2013), 1--23.

\bibitem[ACP12]{acp}
D.~Abramovich, L.~Caporaso, and S.~Payne, \emph{The tropicalization of the
  moduli space of curves}, preprint arXiv:1212.0373, 2012.

\bibitem[AM69]{AtiyahMacdonald69}
M.~Atiyah and I.~Macdonald, \emph{Introduction to commutative algebra},
  Addison-Wesley Publishing Co., Reading, Mass.-London-Don Mills, Ont., 1969.

\bibitem[Bak08a]{Baker08b}
M.~Baker, \emph{An introduction to {B}erkovich analytic spaces and
  non-{A}rchimedean potential theory on curves}, {$p$}-adic geometry, Univ.
  Lecture Ser., vol.~45, Amer. Math. Soc., Providence, RI, 2008, pp.~123--174.

\bibitem[Bak08b]{Baker08}
\bysame, \emph{Specialization of linear systems from curves to graphs}, Algebra
  Number Theory \textbf{2} (2008), no.~6, 613--653.

\bibitem[Ber90]{Berkovich90}
V.~Berkovich, \emph{Spectral theory and analytic geometry over
  non-{A}rchimedean fields}, Mathematical Surveys and Monographs, vol.~33,
  American Mathematical Society, Providence, RI, 1990.

\bibitem[Ber93]{Berkovich93}
\bysame, \emph{\'{E}tale cohomology for non-{A}rchimedean analytic spaces},
  Inst. Hautes \'Etudes Sci. Publ. Math. (1993), no.~78, 5--161 (1994).

\bibitem[Ber99]{Berkovich99}
\bysame, \emph{Smooth {$p$}-adic analytic spaces are locally contractible},
  Invent. Math. \textbf{137} (1999), no.~1, 1--84.

\bibitem[Ber00]{Berkovich00}
\bysame, \emph{An analog of {T}ate's conjecture over local and finitely
  generated fields}, Internat. Math. Res. Notices (2000), no.~13, 665--680.

\bibitem[Ber04]{Berkovich04}
\bysame, \emph{Smooth {$p$}-adic analytic spaces are locally contractible.
  {II}}, Geometric aspects of Dwork theory. Vol. I, II, Walter de Gruyter GmbH
  \& Co. KG, Berlin, 2004, pp.~293--370.

\bibitem[BGR84]{BGR84}
S.~Bosch, U.~G{\"u}ntzer, and R.~Remmert, \emph{Non-{A}rchimedean analysis},
  Grundlehren der Mathematischen Wissenschaften, vol. 261, Springer-Verlag,
  Berlin, 1984.

\bibitem[BL93a]{BoschLutkebohmert93}
S.~Bosch and W.~L\"{u}tkebohmert, \emph{Formal and rigid geometry. {I}. {R}igid
  spaces}, Math. Ann. \textbf{295} (1993), no.~2, 291--317.

\bibitem[BL93b]{BoschLutkebohmert93b}
\bysame, \emph{Formal and rigid geometry. {II}. {F}lattening techniques}, Math.
  Ann. \textbf{296} (1993), no.~3, 403--429.

\bibitem[BLR95a]{BLR95}
S.~Bosch, W.~L\"{u}tkebohmert, and M.~Raynaud, \emph{Formal and rigid geometry.
  {III}. {T}he relative maximum principle}, Math. Ann. \textbf{302} (1995),
  no.~1, 1--29.

\bibitem[BLR95b]{BLR95b}
\bysame, \emph{Formal and rigid geometry. {IV}. {T}he reduced fibre theorem},
  Invent. Math. \textbf{119} (1995), no.~2, 361--398.

\bibitem[BN07]{BakerNorine07}
M.~Baker and S.~Norine, \emph{{R}iemann-{R}och and {A}bel-{J}acobi theory on a
  finite graph}, Adv. Math. \textbf{215} (2007), no.~2, 766--788.

\bibitem[BN09]{BakerNorine09}
\bysame, \emph{Harmonic morphisms and hyperelliptic graphs}, Int. Math. Res.
  Not. (2009), no.~15, 2914--2955.

\bibitem[Bot59]{Bott59}
R.~Bott, \emph{On a theorem of {L}efschetz}, Michigan Math. J. \textbf{6}
  (1959), 211--216.

\bibitem[BPR11]{BPR11}
M.~Baker, S.~Payne, and J.~Rabinoff, \emph{Nonarchimedean geometry,
  tropicalization, and metrics on curves}, preprint, arXiv:1104.0320v1, 2011.

\bibitem[BR13]{BakerRabinoff13}
M.~Baker and J.~Rabinoff, \emph{The skeleton of the {J}acobian, the {J}acobian
  of the skeleton, and lifting meromorphic functions from tropical to algebraic
  curves}, To appear in Int. Math. Res. Not. arXiv:1308.3864, 2013.

\bibitem[Cap12]{Caporaso12}
L.~Caporaso, \emph{Gonality of algebraic curves and graphs}, preprint
  arXiv:1201.6246v3, 2012.

\bibitem[CC12]{CastryckCools12}
W.~Castryck and F.~Cools, \emph{Newton polygons and curve gonalities}, J.
  Algebraic Combin. \textbf{35} (2012), no.~3, 345--366.

\bibitem[CDPR12]{tropicalBN}
F.~Cools, J.~Draisma, S.~Payne, and E.~Robeva, \emph{A tropical proof of the
  {B}rill-{N}oether theorem}, Adv. Math. \textbf{230} (2012), no.~2, 759--776.

\bibitem[CGP14]{cgp}
M.~Chan, S.~Galatius, and S.~Payne, in preparation, 2014.

\bibitem[CKK12]{CornelissenKatoKool12}
G.~Cornelissen, F.~Kato, and J.~Kool, \emph{A combinatorial {L}i-{Y}au
  inequality and rational points on curves}, preprint arXiv:1211.2681, 2012.

\bibitem[Dan75]{Danilov75}
V.~Danilov, \emph{Polyhedra of schemes and algebraic varieties}, Mat. Sb.
  (N.S.) \textbf{139} (1975), no.~1, 146--158, 160.

\bibitem[Del80]{Deligne80}
P.~Deligne, \emph{La conjecture de {W}eil. {II}}, Inst. Hautes \'Etudes Sci.
  Publ. Math. (1980), no.~52, 137--252.

\bibitem[dFKX13]{deFernexKollarXu13}
T.~de~Fernex, J.~Koll\'{a}r, and C.~Xu, \emph{The dual complex of
  singularities}, To appear in proceedings of the conference in honor of Yujiro
  Kawamata's 60th birthday, Adv. Stud. Pure Math. arXiv:1212.1675v2, 2013.

\bibitem[dJ96]{deJong96}
A.~J. de~Jong, \emph{Smoothness, semi-stability and alterations}, Inst. Hautes
  \'Etudes Sci. Publ. Math. (1996), no.~83, 51--93.

\bibitem[DL98]{DenefLoeser98}
J.~Denef and F.~Loeser, \emph{Motivic {I}gusa zeta functions}, J. Algebraic
  Geom. \textbf{7} (1998), no.~3, 505--537.

\bibitem[Duc13]{Ducros13}
A.~Ducros, \emph{Les espaces de {B}erkovich sont mod\'er\'es (d'apr\`es {E}hud
  {H}rushovski et {F}ran\c cois {L}oeser)}, Ast\'erisque (2013), no.~352, Exp.
  No. 1056, x, 459--507, S{\'e}minaire Bourbaki. Vol. 2011/2012. Expos{\'e}s
  1043--1058.

\bibitem[FGP12]{limits}
T.~Foster, P.~Gross, and S.~Payne, \emph{Limits of tropicalizations}, To appear
  in Israel J. Math. arXiv:1211.2718, 2012.

\bibitem[GK08]{GathmannKerber08}
A.~Gathmann and M.~Kerber, \emph{A {R}iemann-{R}och theorem in tropical
  geometry}, Math. Z. \textbf{259} (2008), no.~1, 217--230.

\bibitem[Gub13]{Gubler13}
W.~Gubler, \emph{A guide to tropicalizations}, Algebraic and combinatorial
  aspects of tropical geometry, Contemp. Math., vol. 589, Amer. Math. Soc.,
  Providence, RI, 2013, pp.~125--189.

\bibitem[Hac08]{Hacking08}
P.~Hacking, \emph{The homology of tropical varieties}, Collect. Math.
  \textbf{59} (2008), no.~3, 263--273.

\bibitem[Har77]{Hartshorne77}
R.~Hartshorne, \emph{Algebraic geometry}, Springer-Verlag, New York, 1977,
  Graduate Texts in Mathematics, No. 52.

\bibitem[HHM08]{HaskellHrushovskiMacpherson08}
D.~Haskell, E.~Hrushovski, and D.~Macpherson, \emph{Stable domination and
  independence in algebraically closed valued fields}, Lecture Notes in Logic,
  vol.~30, Association for Symbolic Logic, Chicago, IL, 2008.

\bibitem[HK06]{HrushovskiKazhdan06}
E.~Hrushovski and D.~Kazhdan, \emph{Integration in valued fields}, Algebraic
  geometry and number theory, Progr. Math., vol. 253, Birkh\"auser Boston,
  Boston, MA, 2006, pp.~261--405.

\bibitem[HK12]{HelmKatz12}
D.~Helm and E.~Katz, \emph{Monodromy filtrations and the topology of tropical
  varieties}, Canad. J. Math. \textbf{64} (2012), no.~4, 845--868.

\bibitem[HL11]{HrushovskiLoeser11}
E.~Hrushovski and F.~Loeser, \emph{Monodromy and the {L}efschetz fixed point
  formula}, preprint, arXiv:1111.1954, 2011.

\bibitem[HL12]{HrushovskiLoeser10}
\bysame, \emph{Nonarchimedean topology and stably dominated types}, preprint,
  arXiv:1009.0252v3, 2012.

\bibitem[HLP12]{HrushovskiLoeserPoonen12}
E.~Hrushovski, F.~Loeser, and B.~Poonen, \emph{Berkovich spaces embed in
  euclidean spaces}, preprint, arXiv:1210.6485, 2012.

\bibitem[Igu75]{Igusa75}
J.-I. Igusa, \emph{Complex powers and asymptotic expansions. {II}. {A}symptotic
  expansions}, J. Reine Angew. Math. \textbf{278/279} (1975), 307--321.

\bibitem[Jon12]{Jonsson12}
M~Jonsson, \emph{Dynamics on {B}erkovich spaces of low dimensions}, To appear
  in {B}erkovich spaces and applications. {S}\'{e}minaires et {C}ongr\`{e}s.
  arXiv:1201.1944v1, 2012.

\bibitem[JP14]{tropicalGP}
D.~Jensen and S.~Payne, \emph{Tropical multiplication maps and the
  {G}ieseker-{P}etri {T}heorem}, preprint, arXiv:1401.2584, 2014.

\bibitem[KK11]{KapovichKollar11}
M.~Kapovich and J.~Koll{\'a}r, \emph{Fundamental groups of links of isolated
  singularities}, To appear in J. Amer. Math. Soc. arXiv:1109.4047v1, 2011.

\bibitem[KKMSD73]{KKMS}
G.~Kempf, F.~Knudsen, D.~Mumford, and B.~Saint-Donat, \emph{Toroidal
  embeddings. {I}}, Lect. Notes in Math., vol. 339, Springer-Verlag, Berlin,
  1973.

\bibitem[Kol13]{Kollar13}
J.~Koll\'{a}r, \emph{Links of complex analytic singularities}, Surveys in
  differential geometry. {G}eometry and topology, Surv. Differ. Geom., vol.~18,
  Int. Press, Somerville, MA, 2013, pp.~157--193.

\bibitem[KS06]{KontsevichSoibelman06}
M.~Kontsevich and Y.~Soibelman, \emph{Affine structures and non-{A}rchimedean
  analytic spaces}, The unity of mathematics, Progr. Math., vol. 244,
  Birkh\"auser Boston, Boston, MA, 2006, pp.~321--385.

\bibitem[KS12a]{KatzStapledon12}
E.~Katz and A.~Stapledon, \emph{Tropical geometry and the motivic nearby
  fiber}, Compos. Math. \textbf{148} (2012), no.~1, 269--294.

\bibitem[KS12b]{KerzSaito12}
M.~Kerz and S.~Saito, \emph{Cohomological {H}asse principle and motivic
  cohomology for arithmetic schemes}, Publ. Math. Inst. Hautes \'Etudes Sci.
  (2012), 123--183.

\bibitem[KT02]{KontsevichTschinkel02}
M.~Kontsevich and Y.~Tschinkel, \emph{Nonarchimedean {K}{\"a}hler geometry},
  unpublished, 2002.

\bibitem[LS03]{LoeserSebag03}
F.~Loeser and J.~Sebag, \emph{Motivic integration on smooth rigid varieties and
  invariants of degenerations}, Duke Math. J. \textbf{119} (2003), no.~2,
  315--344.

\bibitem[Mar02]{Marker02}
D.~Marker, \emph{Model theory}, Graduate Texts in Mathematics, vol. 217,
  Springer-Verlag, New York, 2002, An introduction.

\bibitem[Mil68]{Milnor68}
J.~Milnor, \emph{Singular points of complex hypersurfaces}, Annals of
  Mathematics Studies, No. 61, Princeton University Press, Princeton, N.J.,
  1968.

\bibitem[MN12]{MustataNicaise12}
M.~Musta\c{t}\u{a} and J.~Nicaise, \emph{Weight functions on non-archimedean
  analytic spaces and the {K}ontsevich-{S}oibelman skeleton}, preprint,
  arXiv:1212.6328, 2012.

\bibitem[MZ08]{MikhalkinZharkov08}
G.~Mikhalkin and I.~Zharkov, \emph{Tropical curves, their {J}acobians and theta
  functions}, Curves and abelian varieties, Contemp. Math., vol. 465, Amer.
  Math. Soc., Providence, RI, 2008, pp.~203--230.

\bibitem[Nic11]{Nicaise11}
J.~Nicaise, \emph{Singular cohomology of the analytic {M}ilnor fiber, and mixed
  {H}odge structure on the nearby cohomology}, J. Algebraic Geom. \textbf{20}
  (2011), no.~2, 199--237.

\bibitem[NS07]{NicaiseSebag07}
J.~Nicaise and J.~Sebag, \emph{Motivic {S}erre invariants, ramification, and
  the analytic {M}ilnor fiber}, Invent. Math. \textbf{168} (2007), no.~1,
  133--173.

\bibitem[NX13]{NicaiseXu13}
J.~Nicaise and C.~Xu, \emph{The essential skeleton of a degeneration of
  algebraic varieties}, arXiv:1307.4041, 2013.

\bibitem[Pay09]{analytification}
S.~Payne, \emph{Analytification is the limit of all tropicalizations}, Math.
  Res. Lett. \textbf{16} (2009), no.~3, 543--556.

\bibitem[Pay13]{boundarycx}
S.~Payne, \emph{Boundary complexes and weight filtrations}, Michigan Math. J.
  \textbf{62} (2013), 293--322.

\bibitem[Ray74]{Raynaud74b}
M.~Raynaud, \emph{G\'eom\'etrie analytique rigide d'apr\`es {T}ate,
  {K}iehl,{$\ldots $}}, Table {R}onde d'{A}nalyse non archim\'edienne ({P}aris,
  1972), Soc. Math. France, Paris, 1974, pp.~319--327. Bull. Soc. Math. France,
  M\'em. No. 39--40.

\bibitem[Smi14]{Smith14}
G.~Smith, \emph{Brill-{N}oether theory of curves in toric surfaces}, preprint
  arXiv:1403.2317. To appear in J. Pure Appl. Alg., 2014.

\bibitem[Ste06]{Stepanov06}
D.~Stepanov, \emph{A remark on the dual complex of a resolution of
  singularities}, Uspekhi Mat. Nauk \textbf{61} (2006), no.~1(367), 185--186.

\bibitem[Thu05]{ThuillierThesis}
A.~Thuillier, \emph{Th{\'e}orie du potentiel sur ler courbes en
  g{\'e}om{\'e}trie analytique non archim{\'e}dienne. {A}pplications {\`a} la
  the{\'e}orie d'{A}rakelov}, Ph.D. thesis, University of Rennes, 2005.

\bibitem[Thu07]{Thuillier07}
\bysame, \emph{G\'eom\'etrie toro\"\i dale et g\'eom\'etrie analytique non
  archim\'edienne}, Manuscripta Math. \textbf{123} (2007), no.~4, 381--451.

\end{thebibliography}
\end{document}